\documentclass[12pt]{amsart}
\usepackage{fullpage,url,amssymb,enumerate,colonequals}
\usepackage[all,pdf]{xy} 

\usepackage{comment}
\usepackage{color}
\usepackage{charter,euler} 
\DeclareSymbolFont{bchoperators}{T1}{bch}{m}{n}
\SetSymbolFont{bchoperators}{bold}{T1}{bch}{b}{n}
\makeatletter
\renewcommand{\operator@font}{\mathgroup\symbchoperators}
\makeatother

\usepackage{titlesec} 
\titleformat{\section}{\normalfont\bfseries\filcenter}{\thesection}{1em}{}

\usepackage[pdftex]{epsfig}

\newcommand{\F}{\mathbb{F}}

\newcommand{\PP}{\mathbb{P}}
\newcommand{\Q}{\mathbb{Q}}

\newcommand{\Z}{\mathbb{Z}}

\newcommand{\calL}{\mathcal{L}}

\newcommand{\Pic}{\operatorname{Pic}}

\numberwithin{equation}{section}

\newtheorem{theorem}{Theorem}

\definecolor{darkgreen}{rgb}{0,0.5,0}
\definecolor{rem}{rgb}{0.8,0,0}
\definecolor{new}{rgb}{0.7,0,0.6}
\definecolor{reply}{rgb}{0,0,0.8}
\usepackage[
        colorlinks, citecolor=darkgreen,
        backref,
        pdfauthor={Peter Bruin, Maarten Derickx, Michael Stoll},
        pdftitle={Elliptic curves with a point of order 13 defined over cyclic cubic fields},
]{hyperref}
\usepackage[alphabetic,backrefs,lite]{amsrefs} 

\begin{document}

\title{Elliptic curves with a point of order 13\\
       defined over cyclic cubic fields}

\author{Peter Bruin}
\address{Universiteit Leiden,
         Mathematisch Instituut,
         Postbus 9512,
         2300 RA Leiden,
         The Netherlands.}
\email{P.J.Bruin@math.leidenuniv.nl}
\urladdr{\url{http://pub.math.leidenuniv.nl/~bruinpj/}}

\author{Maarten Derickx}
\address{Universiteit Leiden,
	Mathematisch Instituut,
	Postbus 9512,
	2300 RA Leiden,
	The Netherlands.}
\email{maarten@mderickx.nl}
\urladdr{\url{http://www.maartenderickx.nl/}}

\author{Michael Stoll}
\address{Mathematisches Institut,
         Universit\"at Bayreuth,
         95440 Bayreuth, Germany.}
\email{Michael.Stoll@uni-bayreuth.de}
\urladdr{\url{http://www.mathe2.uni-bayreuth.de/stoll/}}

\date{January 14, 2021}

\begin{abstract}
  We show that there is essentially a unique elliptic curve~$E$ defined
  over a cubic Galois extension~$K$ of~$\Q$ with a $K$-rational point
  of order~$13$ and such that $E$ is not defined over~$\Q$.
\end{abstract}


\maketitle

\vspace{-1.2em}

\section{Introduction}

Let $K$ be a number field and let $E$ be an elliptic curve over~$K$.
The question what the possible orders of torsion points in $E(K)$ are
(or, more generally, which finite abelian groups occur
as the group of $K$-rational torsion points of some~$E$) has received
much attention in the recent past. Mazur~\cite{Mazur1977} famously
solved this problem for $K = \Q$. Kenku and Momose~\cite{KenkuMomose}
and Kamienny~\cites{Kamienny1986,Kamienny1992} dealt with quadratic number fields.
Merel~\cite{Merel1996} proved that for fields~$K$ of given degree~$d$,
there are only finitely many possibilities for the order of a torsion
point (and therefore also for the torsion subgroup). Jeon, with various
coauthors~\cites{JeonKimSchweizer,JeonKimLee-cubic-families}
determined which torsion structures occur infinitely often for
cubic fields.
Najman~\cite{Najman2016} found a sporadic example with a point
of order~$21$ over a cubic field (with the curve defined over~$\Q$);
this is not on Jeon et al.'s list. The second author together with Etropolski, van Hoeij, Morrow and Zureick-Brown~\cites{DEHMZB}, building on the results of Parent~\cites{Parent2000,Parent2003}, proved that $\Z/21\Z$ is actually the only torsion structure that occurs finitely often over cubic fields. Jeon~\cite{Jeon2016} also determined
which torsion structures occur infinitely often over cyclic cubic fields.
The second author and Najman~\cite{DerickxNajman} classified all torsion
groups that occur over cyclic cubic fields.
There are similar results by Jeon and coauthors for quartic
fields~\cites{JeonKimPark,JeonKimLee-quartic,JeonKimLee-dihedral1,JeonKimLee-dihedral2}.
Results on which prime numbers can occur as the order of a torsion point
over a field of degree $d \le 7$ can be found in a forthcoming
paper by Kamienny, Stein and the last two authors of this note~\cite{DKSS}.

In this note, we consider the cubic case. More precisely, we complete
the classification of elliptic curves over cyclic cubic fields
that have a point of order~$13$. Jeon in~\cite{Jeon2016} already found
an infinite family (with parameter space an open subset of the projective line)
of such curves. They are obtained by pulling back rational points under
a cyclic degree~$3$ Galois cover $X_1(13) \to \PP^1$, which is derived
from the action of the diamond operators on~$X_1(13)$. This implies that
the target~$\PP^1$ is a modular curve itself, and
all the curves in the family are in fact already defined over~$\Q$
(and acquire a point of order~$13$ over a cyclic cubic extension).
We show that outside this family, there is essentially one other example,
which is an elliptic curve that cannot be defined over~$\Q$.
This can be seen as a complement to~\cite{DerickxNajman}, where similar
results are obtained for points of order $16$ and~$20$.

\subsection*{Acknowledgments} \strut

We thank the organizers of the conference on ``Torsion groups and Galois
representations of elliptic curves'' in Zagreb in June~2018, where the work
described in this note was done, and Daeyeol Jeon for providing the
motivation for this work by pointing out in his talk at the conference
that $X_1(13)$ arises as a Galois covering of~$\PP^1$ of degree~$3$
over~$\Q$, raising the question whether there are additional ``sporadic''
elliptic curves over a cyclic number field with a point of order~$13$.

The computations were done using the Magma computer algebra system~\cite{Magma}.


\section{The result}

Our goal is to classify all elliptic curves~$E$ defined over a cyclic
cubic extension~$K$ of~$\Q$ such that $E(K)$ contains points of order~$13$.
The main result is as follows.

\begin{theorem} \label{T:main}
  Let $K$ be a cubic Galois extension of~$\Q$ and let $E$ be an elliptic
  curve defined over~$K$ with $E(K)[13] \neq 0$. Then either $E$ is defined
  over~$\Q$, or else $K = \Q(\alpha)$ with
  \[ \alpha^3 - \alpha^2 - 82 \alpha + 64 = 0 \]
  and $E$ is isomorphic to a Galois conjugate of the curve
  \[ E_0 \colon y^2 + (1-c) xy - b y = x^3 - b x^2 \,, \]
  where
  \[ b = \frac{10 \alpha^2 + 90 \alpha - 1936}{19773} \quad\text{and}\quad
     c = \frac{6 \alpha^2 + 50 \alpha - 208}{1521} \,.
  \]
\end{theorem}

To obtain this result, we find all degree~$3$ morphisms $X_1(13) \to \PP^1$
that are defined over~$\Q$ and determine all their fibers above rational
points that give rise to a cyclic cubic extension of~$\Q$. There are
exactly~$13$ such morphisms (up to automorphisms of~$\PP^1$). One of these
is obtained by dividing by a subgroup of order~$3$ of the group generated
by the diamond operators on~$X_1(13)$. All its fibers are cyclic or split;
this gives rise to the family of elliptic curves over~$\Q$ with points of
order~$13$ defined over a cyclic cubic field found by Jeon~\cite{Jeon2016}.
Explicitly, this family can be obtained as
\[ E_t \colon y^2 = x^3 - 27 A(t) x + 54 (t^2 + 1) B(t) \]
with
\begin{align*}
  A(t) &= \frac{t^8 - 5 t^7 + 7 t^6 - 5 t^5 + 5 t^3 + 7 t^2 + 5 t + 1}%
               {t^4 - t^3 + 5 t^2 + t + 1} \qquad\text{and} \\
  B(t) &= \frac{t^{12} - 8 t^{11} + 25 t^{10} - 44 t^9 + 40 t^8 + 18 t^7 - 40 t^6 - 18 t^5
                       + 40 t^4 + 44 t^3 + 25 t^2 + 8 t + 1}%
               {(t^4 - t^3 + 5 t^2 + t + 1)^2} \,.
\end{align*}
A point of order~$13$ on~$E_t$ is given by
\[  P_t = \Bigl(\frac{36 t w + 3 (t^6 - 3 t^5 + 4 t^4 - 6 t^3 - 8 t^2 + 3 t + 1)}%
                     {t^4 - t^3 + 5 t^2 + t + 1},
                \frac{108 t ((t + 1) w - t)}{t^4 - t^3 + 5 t^2 + t + 1} \Bigr) \,,
\]
where
\[ w^3 + (-t^3 + t^2 - 3 t + 1) w^2 + (-t^3 + 2 t^2 - 2 t) w + t^2 = 0 \,. \]
This last polynomial has discriminant~$t^4 (t^4 - t^3 + 5 t^2 + t + 1)^2$
and therefore defines, for $t \in \Q \setminus \{0\}$, a cyclic cubic
number field. (Our parameter~$t$ is related to $t_\text{Jeon}$ of~\cite{Jeon2016}
by $t_\text{Jeon} = -\frac{7}{72} - \frac{1}{36 t}$.)

The other~$12$ morphisms
fall into two orbits under the group of diamond operators (which is cyclic
of order~$6$). These morphisms are not Galois coverings of~$\PP^1$, so
their fibers usually define $S_3$-extensions of~$\Q$. The condition for
a fiber to be cyclic is expressed by requiring the discriminant of a cubic
polynomial over~$\Q(t)$ (where $t$ is a parameter on~$\PP^1$; adjoining
a root of the cubic defines the covering) to be a square. This defines a
hyperelliptic curve. For one of the two orbits, we obtain a curve of genus~$2$,
for which we can prove using Chabauty's method that it has exactly five
rational points. Three of these points arise from split fibers containing
cusps, but one pair of points corresponds to a rational point on~$\PP^1$
that has a cyclic fiber above it, leading to the curve~$E_0$ in Theorem~\ref{T:main}.
The remaining orbit leads to a curve of genus~$3$, for which we can prove
that it has exactly three rational points; they all correspond to
ramified fibers containing cusps. The details are given in the next section.


\section{Proof of the theorem}

In the following, we will write $X$ for~$X_1(13)$.
We use the model of~$X$ given by
\[ y^2 + (x^3 + x^2 + 1) y = x^2 + x \,; \]
see Sutherland's table~\cite{SutherlandTable}.
In particular, $X$ has genus~$2$ and is (therefore) hyperelliptic.
The \emph{canonical class} on~$X$ is the linear equivalence class
of divisors arising by pulling back a point under the hyperelliptic
covering map $X \to \PP^1$; it is the same as the class containing
the canonical divisors.

The map $X = X_1(13) \to X_0(13)$ obtained by sending a point of order~$13$
to the group it generates is Galois over~$\Q$.
Its automorphism group consists of the diamond operators and is
canonically isomorphic to $(\Z/13\Z)^\times/\{\pm 1\}$, which is a cyclic
group of order~$6$; we denote this group by~$G$.

It is known that $X(\Q)$ consists of the six rational cusps
and that the group of rational points on its Jacobian~$J$
is cyclic of order~$19$; see~\cite{MazurTate}. The rational points on~$X$
form one orbit under~$G$;
they are the two points at infinity and the points $(-1,-1)$, $(-1,0)$, $(0,-1)$
and~$(0,0)$ on our model of~$X$.

Since the genus of~$X$ is~$2$,
every rational point on~$J$ except the origin has a unique representation
as an effective divisor~$D$ of degree~$2$ minus the canonical class,
with the points in the support of~$D$ either rational or defined over a
quadratic extension of~$\Q$ and conjugate, where $D$ is not in the canonical class.
Since the six rational points lead to exactly~$18$ effective divisors of degree~$2$
outside the canonical class, they account for all the rational points in~$J$,
which implies that there are no quadratic points with irrational $x$-coordinate.
See~\cites{BBDN,DMK}.

Now consider a point~$P \in X$ with $[\Q(P) : \Q] = 3$. The sum of~$P$
and its two Galois conjugates is a rational effective divisor~$D$ of degree~$3$,
so it gives a rational point on the symmetric cube~$S$ of~$X$. Under the
canonical map $\pi \colon S \to \Pic^3(X)$, it maps to a rational point on~$\Pic^3(X)$.
The fibers of~$\pi$ are $\PP^1$'s; this follows from the Riemann-Roch theorem.

Since $\Pic^3(X)$ is isomorphic to~$J$ (note that $X$ has rational points),
it has exactly~$19$ rational points. Therefore $D$ lies in the fiber of~$\pi$ above
one of these points. Six of the rational points on~$\Pic^3(X)$ arise as
a rational point on~$X$ plus the canonical class. This implies that all
divisors in the corresponding fiber contain this rational point and can therefore
never contain a cubic point.

For the remaining $13$ rational points on~$\Pic^3(X)$, the corresponding
line bundle~$\calL$ is base\-point-free. In this case, the fiber of~$\pi$ above
the point can be identified with the target of the morphism $X \to \PP^1$ defined by
the two-dimensional space of global sections of the line bundle~$\calL$.
This implies that each cubic point on~$X$ lies in the fiber of one of
these morphisms $X \to \PP^1$ above a rational point of~$\PP^1$.

These $13$ rational points on~$\Pic^3(X)$
consist of one point that is fixed by~$G$ and two orbits of size~$6$ under~$G$.
It clearly suffices to determine the cyclic cubic points in the fibers of
the degree~$3$ morphisms $X \to \PP^1$ associated to one representative
of each orbit.

Let $G'$ be the subgroup of order~$3$ of~$G$. The quotient $X/G'$ is a curve
of genus~$0$, so the map $X \to X/G'$ must show up in our list. Since this
quotient is unique, the map must correspond to the point in~$\Pic^3(X)(\Q)$
that is fixed by~$G$. Since $X \to X/G'$ is a Galois covering, all its fibers
over rational points are either ramified, split, or cyclic. The modular curve~$X/G'$
is a double cover of~$X_0(13)$; viewing $G'$ as the group of diamond
operators $\{\langle 1\rangle, \langle 3\rangle, \langle 9\rangle\}$,
we see that $X/G'$ is a fine moduli space outside the branch locus of
$X\to X/G'$ (the image in $X/G'$ of the zero locus of $x^2+x+1$ in our
model of~$X$.) Since the cubic points arising in fibers of $X \to
X/G'$ map to rational points on~$X/G'$ outside this branch locus,
the elliptic curves they represent are defined over~$\Q$. This accounts for
the first alternative in Theorem~\ref{T:main}. (We have made this one-parameter
family explicit in the previous section.)

One representative of one of the other two orbits of degree~$3$ morphisms
to~$\PP^1$ is given by the $y$-coordinate map of our model of~$X$.
The discriminant with respect to~$x$ of the equation defining~$X$ is
\[ d_1(y) = (y + 1) (-27 y^5 - 31 y^4 - 6 y^3 + 6 y^2 + 5 y + 1) \,; \]
the condition that this is a square then defines a hyperelliptic curve~$D_1$
of genus~$2$. A quick search finds five rational points on~$D_1$:
one point with $y = -1$ and two each with $y = 0$ and $y = -\frac{4}{13}$.
So there are three fibers of the $y$-coordinate map with Galois group
contained in $A_3$. The first two contain rational points on~$X$, but the
fiber above~$-\frac{4}{13}$ really consists of three conjugate points
defined over the cyclic extension~$K$; they and the other points in their
$G$-orbits give rise to the curve~$E_0$ and its Galois conjugates mentioned
in~\ref{T:main}. (It can be easily checked that the point~$(0,0)$ on~$E_0$
indeed has order~$13$. The discriminant of~$K$ is $(13 \cdot 19)^2$.)

Using the Magma implementation of $2$-descent on hyperelliptic Jacobians
as described in~\cite{Stoll2001}, we find that the Mordell-Weil rank of the
Jacobian of~$D_1$ is at most~$1$. From the rational points we have found
on~$D_1$, we can easily construct a rational point of infinite order on the Jacobian.
A combination of Chabauty's method with the Mordell-Weil sieve as
explained in~\cite{BruinStoll2010} and implemented in Magma then quickly
proves that the five points we found are indeed all the rational points on~$D_1$.

A representative of the remaining orbit is $X \to \PP^1$ given by~$\frac{y+1}{x}$.
Writing $t$ for the parameter on~$\PP^1$, we have $y = xt - 1$. Plugging
this into the equation of~$X$ and taking the discriminant with respect to~$x$
gives
\[ d_2(t) = t (t+1)^3 (-4 t^5 + 5 t^4 - t^3 - 25 t^2 - 23 t - 4) \,. \]
Setting $d_2(t)/(t+1)^2$ equal to a square gives a hyperelliptic curve~$D_2$
of genus~$3$. It has three obvious rational Weierstrass points at infinity
and with $t = -1$ or~$0$. We do not find any other rational point.
Using $2$-descent again, we can show that the Mordell-Weil rank of the
Jacobian is~$0$. A minimal model of~$D_2$ is
\[ v^2 + (u^3 + u^2) v = u^7 - 8 u^5 - 13 u^4 - 7 u^3 - 2 u^2 - u \,; \]
from this we see that $D_2$ has good reduction at~$2$. The reduction
has exactly three $\F_2$-points, which are the images of the three
rational points we had found. Since their residue disks are fixed by the
hyperelliptic involution, we know that each of the three residue disks
contains an odd number of rational points. Since the Mordell-Weil group
is finite, all $2$-adic integrals $\int_P^Q \omega$ between rational points
$P, Q \in D_2(\Q)$ for regular differentials~$\omega$ must vanish.
In particular, we can, for each residue class, choose a differential
whose reduction mod~$2$ does not vanish in the corresponding $\F_2$-point
on the reduction. By~\cite{Stoll2006a}*{Section~6}, the corresponding
integral vanishes for at most two points in the residue class. Since it
has to vanish at each rational point and the number of rational points
in the residue class is odd, there is only the known rational point
in each of the three residue classes, which shows that $\#D_2(\Q) = 3$.
These three points are all images of cusps, so we do not obtain any
further cyclic cubic points on~$X$. This concludes the proof.



\begin{bibdiv}
\begin{biblist}

\bib{Magma}{article}{
   author={Bosma, Wieb},
   author={Cannon, John},
   author={Playoust, Catherine},
   title={The Magma algebra system. I. The user language},
   note={Computational algebra and number theory (London, 1993)},
   journal={J. Symbolic Comput.},
   volume={24},
   date={1997},
   number={3-4},
   pages={235--265},
   issn={0747-7171},
   review={\MR{1484478}},
   doi={10.1006/jsco.1996.0125},
}

\bib{BBDN}{article}{
   author={Bosman, Johan},
   author={Bruin, Peter},
   author={Dujella, Andrej},
   author={Najman, Filip},
   title={Ranks of elliptic curves with prescribed torsion over number
   fields},
   journal={Int. Math. Res. Not. IMRN},
   date={2014},
   number={11},
   pages={2885--2923},
   issn={1073-7928},
   review={\MR{3214308}},
   doi={10.1093/imrn/rnt013},
}

\bib{BruinStoll2010}{article}{
   author={Bruin, Nils},
   author={Stoll, Michael},
   title={The Mordell-Weil sieve: proving non-existence of rational points on curves},
   journal={LMS J. Comput. Math.},
   volume={13},
   date={2010},
   pages={272--306},
   issn={1461-1570},
   review={\MR{2685127}},
   doi={10.1112/S1461157009000187},
}


\bib{DEHMZB}{misc}{
	author={Derickx, Maarten},
	author={Etropolski, Anastassia},
	author={van Hoeij, Mark},
	author={Morrow, Jackson S.},
	author={Zureick-Brown, David},
	title={Sporadic Cubic Torsion},
	date={2020-07-28},
	note={Preprint, \texttt{arXiv:2007.13929}},
	eprint={arXiv:2007.13929},
}

\bib{DerickxNajman}{article}{
	author={Derickx, Maarten},
	author={Najman, Filip},
	title={Torsion of elliptic curves over cyclic cubic fields},
	journal={Math. Comp.},
	volume={88},
	date={2019},
	number={319},
	pages={2443--2459},
	issn={0025-5718},
	review={\MR{3957900}},
	doi={10.1090/mcom/3408},
}

\bib{DKSS}{misc}{
   author={Derickx, Maarten},
   author={Kamienny, Sheldon},
   author={Stein, William},
   author={Stoll, Michael},
   title={Torsion points on elliptic curves over number fields of small degree},
   date={2017-07-02},
   note={Preprint, \texttt{arXiv:1707.00364}},
   eprint={arXiv:1707.00364},
}


\bib{DMK}{article}{
   author={Derickx, Maarten},
   author={Mazur, Barry},
   author={Kamienny, Sheldon},
   title={Rational families of 17-torsion points of elliptic curves over
   number fields},
   conference={
      title={Number theory related to modular curves---Momose memorial
      volume},
   },
   book={
      series={Contemp. Math.},
      volume={701},
      publisher={Amer. Math. Soc., Providence, RI},
   },
   date={2018},
   pages={81--104},
   review={\MR{3755909}},
   doi={10.1090/conm/701/14142},
}

\bib{Jeon2016}{article}{
   author={Jeon, Daeyeol},
   title={Families of elliptic curves over cyclic cubic number fields with
   prescribed torsion},
   journal={Math. Comp.},
   volume={85},
   date={2016},
   number={299},
   pages={1485--1502},
   issn={0025-5718},
   review={\MR{3454372}},
   doi={10.1090/mcom/3012},
}

\bib{JeonKimLee-cubic-families}{article}{
   author={Jeon, Daeyeol},
   author={Kim, Chang Heon},
   author={Lee, Yoonjin},
   title={Families of elliptic curves over cubic number fields with
   prescribed torsion subgroups},
   journal={Math. Comp.},
   volume={80},
   date={2011},
   number={273},
   pages={579--591},
   issn={0025-5718},
   review={\MR{2728995}},
   doi={10.1090/S0025-5718-10-02369-0},
}

\bib{JeonKimLee-quartic}{article}{
   author={Jeon, Daeyeol},
   author={Kim, Chang Heon},
   author={Lee, Yoonjin},
   title={Families of elliptic curves over quartic number fields with
   prescribed torsion subgroups},
   journal={Math. Comp.},
   volume={80},
   date={2011},
   number={276},
   pages={2395--2410},
   issn={0025-5718},
   review={\MR{2813367}},
   doi={10.1090/S0025-5718-2011-02493-2},
}

\bib{JeonKimLee-dihedral1}{article}{
   author={Jeon, Daeyeol},
   author={Kim, Chang Heon},
   author={Lee, Yoonjin},
   title={Infinite families of elliptic curves over dihedral quartic number
   fields},
   journal={J. Number Theory},
   volume={133},
   date={2013},
   number={1},
   pages={115--122},
   issn={0022-314X},
   review={\MR{2981403}},
   doi={10.1016/j.jnt.2012.06.014},
}

\bib{JeonKimLee-dihedral2}{article}{
   author={Jeon, Daeyeol},
   author={Kim, Chang Heon},
   author={Lee, Yoonjin},
   title={Families of elliptic curves with prescribed torsion subgroups over
   dihedral quartic fields},
   journal={J. Number Theory},
   volume={147},
   date={2015},
   pages={342--363},
   issn={0022-314X},
   review={\MR{3276329}},
   doi={10.1016/j.jnt.2014.07.014},
}

\bib{JeonKimPark}{article}{
   author={Jeon, Daeyeol},
   author={Kim, Chang Heon},
   author={Park, Euisung},
   title={On the torsion of elliptic curves over quartic number fields},
   journal={J. London Math. Soc. (2)},
   volume={74},
   date={2006},
   number={1},
   pages={1--12},
   issn={0024-6107},
   review={\MR{2254548}},
   doi={10.1112/S0024610706022940},
}

\bib{JeonKimSchweizer}{article}{
   author={Jeon, Daeyeol},
   author={Kim, Chang Heon},
   author={Schweizer, Andreas},
   title={On the torsion of elliptic curves over cubic number fields},
   journal={Acta Arith.},
   volume={113},
   date={2004},
   number={3},
   pages={291--301},
   issn={0065-1036},
   review={\MR{2069117}},
   doi={10.4064/aa113-3-6},
}

\bib{Kamienny1986}{article}{
   author={Kamienny, S.},
   title={Torsion points on elliptic curves over all quadratic fields},
   journal={Duke Math. J.},
   volume={53},
   date={1986},
   number={1},
   pages={157--162},
   issn={0012-7094},
   review={\MR{835802}},
   doi={10.1215/S0012-7094-86-05310-X},
}

\bib{Kamienny1992}{article}{
   author={Kamienny, S.},
   title={Torsion points on elliptic curves and $q$-coefficients of modular
   forms},
   journal={Invent. Math.},
   volume={109},
   date={1992},
   number={2},
   pages={221--229},
   issn={0020-9910},
   review={\MR{1172689}},
   doi={10.1007/BF01232025},
}

\bib{KenkuMomose}{article}{
   author={Kenku, M. A.},
   author={Momose, F.},
   title={Torsion points on elliptic curves defined over quadratic fields},
   journal={Nagoya Math. J.},
   volume={109},
   date={1988},
   pages={125--149},
   issn={0027-7630},
   review={\MR{931956}},
   doi={10.1017/S0027763000002816},
}

\bib{Mazur1977}{article}{
   author={Mazur, B.},
   title={Modular curves and the Eisenstein ideal},
   journal={Inst. Hautes \'Etudes Sci. Publ. Math.},
   number={47},
   date={1977},
   pages={33--186 (1978)},
   issn={0073-8301},
   review={\MR{488287}},
}

\bib{MazurTate}{article}{
   author={Mazur, B.},
   author={Tate, J.},
   title={Points of order $13$ on elliptic curves},
   journal={Invent. Math.},
   volume={22},
   date={1973/74},
   pages={41--49},
   issn={0020-9910},
   review={\MR{0347826}},
   doi={10.1007/BF01425572},
}

\bib{Merel1996}{article}{
   author={Merel, Lo\"\i c},
   title={Bornes pour la torsion des courbes elliptiques sur les corps de
   nombres},
   language={French},
   journal={Invent. Math.},
   volume={124},
   date={1996},
   number={1-3},
   pages={437--449},
   issn={0020-9910},
   review={\MR{1369424}},
   doi={10.1007/s002220050059},
}

\bib{Najman2016}{article}{
   author={Najman, Filip},
   title={Torsion of rational elliptic curves over cubic fields and sporadic
   points on $X_1(n)$},
   journal={Math. Res. Lett.},
   volume={23},
   date={2016},
   number={1},
   pages={245--272},
   issn={1073-2780},
   review={\MR{3512885}},
   doi={10.4310/MRL.2016.v23.n1.a12},
}

\bib{Parent2000}{article}{
   author={Parent, Pierre},
   title={Torsion des courbes elliptiques sur les corps cubiques},
   language={French, with English and French summaries},
   journal={Ann. Inst. Fourier (Grenoble)},
   volume={50},
   date={2000},
   number={3},
   pages={723--749},
   issn={0373-0956},
   review={\MR{1779891}},
}

\bib{Parent2003}{article}{
   author={Parent, Pierre},
   title={No 17-torsion on elliptic curves over cubic number fields},
   language={English, with English and French summaries},
   journal={J. Th\'eor. Nombres Bordeaux},
   volume={15},
   date={2003},
   number={3},
   pages={831--838},
   issn={1246-7405},
   review={\MR{2142238}},
}

\bib{Stoll2001}{article}{
   author={Stoll, Michael},
   title={Implementing 2-descent for Jacobians of hyperelliptic curves},
   journal={Acta Arith.},
   volume={98},
   date={2001},
   number={3},
   pages={245--277},
   issn={0065-1036},
   review={\MR{1829626}},
   doi={10.4064/aa98-3-4},
}

\bib{Stoll2006a}{article}{
   author={Stoll, Michael},
   title={Independence of rational points on twists of a given curve},
   journal={Compos. Math.},
   volume={142},
   date={2006},
   number={5},
   pages={1201--1214},
   issn={0010-437X},
   review={\MR{2264661}},
   doi={10.1112/S0010437X06002168},
}

\bib{SutherlandTable}{misc}{
   author={Sutherland, Andrew},
   title={Optimized equations for $X_1(N)$},
   note={\url{http://math.mit.edu/~drew/X1_optcurves.html}},
}

\end{biblist}
\end{bibdiv}
\end{document}